\documentclass{article}
\usepackage[style=numeric,
    url=false,
    doi=false,
    isbn=false,
    giveninits=true,
    eprint=true]{biblatex}
\AtEveryBibitem{
  \clearfield{note}
  \clearlist{language}
  \clearfield{month}
}
\usepackage[margin=1in]{geometry}
\usepackage[T1]{fontenc}
\usepackage{graphicx}
\usepackage{booktabs}
\usepackage{bibentry}
\usepackage{hyperref}
\usepackage{bm} 
\usepackage{tabularx}
\usepackage{multirow} 
\usepackage{longtable}
\usepackage{amsmath}
\usepackage{caption}
\usepackage{float}
\usepackage{enumitem}
\usepackage{romannum}
\usepackage{adjustbox}
\usepackage{amsthm}
\usepackage{subcaption}
\usepackage{makecell}
\usepackage{amssymb}
\usepackage{mathtools}
\newtheorem{definition}{Definition}
\newtheorem{theorem}{Theorem}

\newtheorem{conjecture}{Conjecture}
\newcounter{example}[section]
\newenvironment{example}[1][]{\refstepcounter{example}\par\medskip
   \noindent \textbf{Example~\theexample. #1} \rmfamily}{\medskip}

\title{Computing Lower Bounds on the Non-negative Rank\\via SAT Solvers}
\author{Tilman Engel, Thomas Kalinowski, Matthias Schymura}
\date{}

\begin{document}
\pagenumbering{arabic}
\maketitle
\begin{abstract}
    Finding the extension complexity (xc) of a polytope is equivalent to finding the non-negative rank of its slack matrix.
    We provide a formulation for the rectangle and refined rectangle covering, bounding the non-negative rank of diverse non-negative matrices.
    While the bounds are known, our formulation enables us to use boolean satisfiability (SAT) solvers, which proves to be a strong tool.
    We obtain improved values for lower bounds on the non-negative rank of multiple matrices. In particular, we determined the xc of some regular polygons.
\end{abstract}
\section{Introduction}

A convex polytope $P$ in $\mathbb{R}^d$ is defined as an intersection of closed half-spaces. This is formalized by a set of inequalities 
\[P\coloneqq\{x:Ax \leq b\},\]
$A\in\mathbb{R}^{m\times d}$, $b\in \mathbb{R}^m$. The intersection of $P$ and a hyperplane $\{x\in \mathbb{R}^d : a_i^T x = b_i\}$, with $a_i$ being a row of $A$, is called facet.
Equivalently, $P$ is the convex hull of a finite set of points $v_1,\dots ,v_n$, called vertices of $P$ (\cite{grunbaum_convex_2003,ziegler_lectures_1995}). The number of facets is called size of $P$. Convex polytopes occur in optimization, describing the set of feasible solutions of a linear program (LP). The cost to solve an LP depends on the size of the underlying polytope. An approach to reduce those costs is to represent the same feasibility set with fewer facets in a higher dimension using a linear projection. We call a higher dimensional polytope $Q$ extended formulation of $P$, if there exists a linear projection $\pi$ with $P = \pi(Q)$. The size of the smallest extended formulation of $P$ is called extension complexity (xc) of $P$. The slack matrix $S\in \mathbb{R}_{\geq0}^{m\times n}$ of $P$ is defined by
\[S_{ij}\coloneqq b_i-a_i^T v_j.\] The most prominent approach to finding the xc of a polytope is to compute the non-negative rank of its slack matrix. Both values are the same, as was established by~\cite{yannakakis_expressing_1991}. Computing the non-negative rank $r_+(S)$ means finding the smallest $r$ so that $S=UV$ with $U$ being an $m \times r$ matrix and $V$ being an $r\times n$ matrix. In general, this task is NP-hard (\cite{vavasis_complexity_2007}). Focusing on special cases is the natural thing to do. Slack matrices of regular $n$-gons are such an example. The rows are circular permutations of one another. However, even for this simple case the xc is only known for few small $n$. When studying extension complexity (xc) of polygons, the rectangular covering is one of the more promising approaches to new bounds (\cite{fiorini_combinatorial_2012}). The non-negative factorisation can be regarded as sum of rank-$1$ matrices. The support of each rank-$1$ matrix is a rectangle. So the number of rectangles needed to cover the support of the slack matrix is at least the number of rank-$1$ matrices needed to factorise $S$. The latter can be higher, because the rectangles view $S$ as a 0--1 matrix and do not take the value of its elements into account. \\
The refined rectangular covering bound (RRCB) is introduced in~\cite{oelze_computing_2014}. It takes the actual values of the elements into account and thus provides a better lower bound. The RRCB for regular $n$-gons has been computed up to $n=13$ in~\cite{vandaele_linear_2017} and for $n=14$ by~\cite{dewez_computational_2022}. For these cases the bound yields the precise values of the xc. Together with their upper bound, this establishes the exact xc for the first time. Computation for larger $n$ was not possible due to time constraints. An alternative to the mixed integer program (MIP) formulation used in~\cite{dewez_computational_2022} are boolean satisfiability (SAT) models and solvers. They use boolean variables and instead of inequalities, the constraints are formulated in conjunctive normal form (CNF). This field made a lot of progress in recent years and has an active competition scene (\cite{satcompetition}). To the best of our knowledge, there have not been any calculations of combinatorial bounds on the non-negative rank using SAT solvers.

\paragraph{Contribution of the Paper}
We compute two lower bounds for a collection of slack matrices, by developing a SAT model out of a standard MIP approach.
With the RRCB we improve the known lower bound of the non-negative rank on several non-negative matrices. For regular $n$-gons with $17\leq n\leq 20$, we determined the exact xc. During the preparation of this manuscript, we became aware of an independent report,
available \href{https://github.com/ajt60gaibb/OpenProblemsInNLA/blob/main/references/holden-nr01-2026-09-13/report.pdf}{here}, which presents the exact xc for theses cases too.
Additionally, the rectangle covering bound (RCB) will be computed, which improves the known lower bound on the icosidodecahedron. 
Our computations support the conjecture made in~\cite{vandaele_linear_2017}, regarding the presented upper bound for regular polygons being exact. The RRCB is the best lower bound for the xc of regular $n$-gons for all computed $n$ and even equals the upper bound, in every case it is known exactly.
\paragraph{Structure of the Paper}
We start by defining the two lower bounds, namely the RCB and the RRCB. We then present a non-linear MIP approach for the RRCB, as it was used by~\cite{baeckelant_computing_2026}. The SAT model is presented as a modification of the MIP model. The RCB models will not be presented separately as they arise by omitting a constraint from the RRCB models. The last section consists of the implementation details and results, closely following the presentation in~\cite{baeckelant_computing_2026} to accommodate easy comparison. We did not rerun their implementations on our hardware, so that a detailed comparison of run times is not very meaningful. Our focus in this study is more on the quantitative level: Which problem sizes are realistically solvable with the respective computational approach?

\section{Preliminaries}

Before describing the models, the refined rectangle covering has to be defined. Let $P$ be a $d$ dimensional polytope, defined by inequalities $Ax \leq b$, $A\in\mathbb{R}^{m\times d}$, $b\in \mathbb{R}^d$. Let $v_j$ be the vertices of $P$ for $j=1,\dots ,n$ and $a_i$ the rows of $A$ for $i=1,\dots ,m$. The slack matrix $S\in \mathbb{R}_{\geq0}^{m\times n}$ of $P$ is defined by \[S_{ij}\coloneqq b_i-a_i^T v_j.\]
When scaling the inequalities so that $||a_i|| = 1$, the element at $(i,j)$ gives the euclidean distance between vertex $v_j$ and the facet with normal vector $a_i$.
    \begin{definition}
        Let $S$ be a slack matrix whose rows and columns are indexed by sets $\mathcal{I}$ and $\mathcal{J}$, respectively.
        The \textbf{support} is defined as $\mathrm{supp}(S) \coloneqq \{(i,j) \in \mathcal{I} \times \mathcal{J} : S_{i,j} > 0\}$.
        \end{definition} 
With the support a rectangle can be defined.
    \begin{definition}
        A set $I\times J$ with $I \subseteq \mathcal{I}$ und $J \subseteq \mathcal{J}$ is called a \textbf{rectangle} of $S$, if $I\times J \subseteq \mathrm{supp}(S)$. 
    \end{definition}
Rectangles are precisely the supports of rank-1 matrices. Because the factorization is non-negative, values from different rank-$1$ matrices cannot cancel each other out under addition.
    \begin{definition}
        A set of $k$ rectangles $R_l$ is called a \textbf{rectangle covering}, if $\mathrm{supp}(S) = \bigcup_{l=1}^k R_l$. The smallest number of rectangles in any rectangle covering of $S$ is called \textbf{rectangle covering number} and denoted by $\mathrm{RCB}(S)$.
    \end{definition}
Any non-negative factorization defines a rectangle covering. But not every rectangle covering of size RCB$(S)$ allows the construction of a same size factorization. The non-zero pattern of $S$ can be covered by RCB$(S)$ rectangles. However, the structure of the values in supp$(S)$ may be too complex for RCB$(S)$ rank-$1$ matrices to capture. This justifies the abbreviation RCB, as it bounds the non-negative rank $r_+$ from below. \\
The idea of a refined rectangle covering was first introduced in~\cite{oelze_computing_2014}.
    \begin{definition}\label{refined}
        A rectangle covering $R_1,\dots,R_k$ of a non-negative matrix $S$ is called a \textbf{refined covering} if \[|\{l:R_l \cap \{(i_1,j_1),(i_2,j_2)\} \neq \emptyset\}| \geq 2\]
        holds for all pairs $(i_1,j_1)$, $(i_2,j_2)$ with $S_{i_1,j_1} \cdot S_{i_2,j_2} > S_{i_1,j_2} \cdot S_{i_2,j_1}$. The smallest size of any refined covering of $S$ is called the \textbf{refined rectangle covering number} and denoted by $\mathrm{RRCB}(S)$.
    \end{definition}
This refinement takes values of $S$ into account, by looking at determinants of all $2\times2$ submatrices. If the non-negative rank is 2, one rank-$1$ matrix is not sufficient to represent its values. A second rectangle is needed, so that the determinant is non-zero.
The RRCB is a tighter bound on the non-negative rank $r_+$ of $S$~\cite{oelze_computing_2014}:
\begin{theorem}
\[\mathrm{RCB}(S) \leq \mathrm{RRCB}(S) \leq r_+(S).\]
\end{theorem}
\begin{example}
An example for which $\mathrm{RCB}(S)< \mathrm{RRCB(S)} < r_+(S)$ is given by \[S =\begin{pmatrix}
    1 & 1 & 1\\
    1 & 2 & 2\\
    1 & 2 & 3
\end{pmatrix}.\]
It has no zeros, so one rectangle can cover the whole matrix, hence $\mathrm{RCB}(S)=1$. The submatrix \[\begin{pmatrix}
    1 & 1\\
    1 & 2
\end{pmatrix}\] has determinant $1$ so it needs two rectangles to cover. Because of the absence of zeros, we can just choose two rectangles, both covering the whole matrix $S$. This way the refined condition is satisfied everywhere. The matrix has rank $3$ so its nonnegative rank must also be $3$. This example illustrates, that with few zeros the RRCB tends to not improve much over the RCB.
\end{example}
\begin{example}
    Another matrix of rank 3 is \[S =\begin{pmatrix}
    1 & 2 & 0\\
    1 & 1 & 1\\
    0 & 1 & 1
\end{pmatrix}.\]
    Here $\mathrm{RCB}(S)< \mathrm{RRCB(S)} \leq r_+(S)$ holds. We need two rectangles to cover the matrix. The refined covering has size three. One rectangle needs to cover element $(2,1)$ and another one $(1,2)$. Both can not cover element $(3,3)$, so we need a third one.
\end{example}

\section{A Non-Linear Integer Program to Find Lower Bounds}\label{sec_mip}

We reduce the problem of determining the RRCB of a given polytope to a sequence of 0--1 feasibility problems. The model has two \textbf{inputs}. There is
    \begin{itemize}
        \item the slack matrix $S \in \mathbb{R}_{\geq 0}^{m\times n}$ of the given polytope and
        \item the number of rectangles $t \in \mathbb{N}$ to cover $S$ with.
    \end{itemize}
If there exists a refined rectangle covering of size $t$, a solution is returned. \\
Otherwise ``unsatisfiable'' is returned. The \textbf{variables} are integers. There are
\begin{itemize}
    \item binary row and column indices of $t$ rectangles $u \in \{0,1\}^{m\times t}$, $v \in \{0,1\}^{n\times t}$ and
    \item auxiliary variables $x \in \{0,1\}^{m\times n\times t}$ with $x_{i,j,k} = u_{i,k} \cdot v_{j,k}$ which indicate, if a given element of $S$ is covered by rectangle $k$.
\end{itemize}
These variables must satisfy three Conditions~\eqref{mip_eq1}-\eqref{mip_eq3} to describe a rectangle covering. For the covering to be refined, we need the Constraint~\eqref{mip_eq4}:
\begin{align}
        x_{i,j,k} = u_{i,k} \cdot v_{j,k} &,\quad \forall (i,j,k) \in [m]\times[n]\times[t],\label{mip_eq1} \\
         \sum_{k=1}^t x_{i,j,k} \geq 1&,\quad \forall (i,j)\in \text{supp}(S) ,\label{mip_eq2} \\
         x_{i,j,k} = 0 &,\quad \forall (i,j)\in [m]\times[n] \setminus \text{supp}(S),\  k\in[t], \label{mip_eq3} \\
        \sum_{k=1}^t (x_{i_1,j_1,k} + x_{i_2,j_2,k} - x_{i_1,j_1,k} \cdot x_{i_2,j_2,k}) \geq 2 &,\quad \forall \{(i_1,j_1),(i_2,j_2)\} \in D. \label{mip_eq4}
\end{align}
We use
\[D \coloneqq \{\{(i_1,j_1),(i_2,j_2)\} \in \left([m]\times [n]\right)^{2}:S_{i_1,j_1} \cdot S_{i_2,j_2} > S_{i_1,j_2} \cdot S_{i_2,j_1}\}\]
as the set of pairs for which the refined condition needs to hold.

The first condition arises from the definition of the auxiliary variables $x$, which makes this program non-linear. This structure models the Cartesian product. Constraint \eqref{mip_eq2} forces every element $(i,j)$ in the support to be covered by at least one rectangle. Condition \eqref{mip_eq3} prevents elements outside the support from being covered. Now to ensure the constraint \[|\{l:R_l \cap \{(i_1,j_1),(i_2,j_2)\} \neq \emptyset\}| \geq 2\] from Definition~\ref{refined} holds, we introduce the Inequalities \eqref{mip_eq4}. Summand $k$ is zero if none of the two elements are covered by rectangle $k$ and one if at least one element is covered by it. This must be true for at least two rectangles.

\section{RRCB SAT Model}

The MIP formulation is already a feasibility problem with binary variables, which we can interpret as boolean. The difficulty in reformulation now lies in the inequalities. Using integers allows counting, as in Condition~\eqref{mip_eq4}. This is not directly possible in a SAT formulation. Our SAT model has the same $2$ \textbf{inputs}:
    \begin{itemize}
        \item the slack matrix $S \in \mathbb{R}_{\geq 0}^{m\times n}$ of the given polytope and
        \item the number of rectangles $t \in \mathbb{N}$ to cover $S$ with.
    \end{itemize}
It also has the same outputs and variables. We stick with the
\begin{itemize}
        \item boolean row and column indices of $t$ rectangles $u \in \{\texttt{F},\texttt{T}\}^{m\times t}$, $v \in \{\texttt{F},\texttt{T}\}^{n\times t}$ and
        \item auxiliary variables $x \in \{\texttt{F},\texttt{T}\}^{m\times n\times t}$ with $x_{i,j,k} = u_{i,k} \land v_{j,k}$.
\end{itemize}
The four types of \textbf{conditions}, the variables must satisfy to be a refined rectangle covering, are formulated below.
    \begin{align}
        x_{i,j,k} \iff u_{i,k} \land v_{j,k} &,\quad \forall (i,j,k) \in [m]\times [n]\times[t]\label{eq1} \\
         \bigvee_{k=1}^t (x_{i,j,k}) &,\quad \forall (i,j)\in \text{supp}(S)\label{eq2} \\
         \lnot(x_{i,j,k}) &,\quad \forall (i,j)\in [m]\times [n] \setminus \text{supp}(S),\  k\in[t] \label{eq3} \\
        \bigwedge_{l=1}^t \bigvee_{\substack{k=1\\ k \neq l}}^t (x_{i_1,j_1,k} \vee x_{i_2,j_2,k}) &,\quad \forall \{(i_1,j_1),(i_2,j_2)\} \in D \label{eq4} 
    \end{align}
Condition \eqref{eq1} arises from the definition of $x$. Every element in the support of $S$ has to be in at least one rectangle, as required by Constraint \eqref{eq2}. If the element is not inside the support, it can't be in any rectangle as enforced in~\eqref{eq3}.\\
For every pair of elements as in Definition~\ref{refined} \[|\{l:R_l \cap \{(i_1,j_1),(i_2,j_2)\} \neq \emptyset\}| \geq 2\] must hold. For $l_1\neq l_2$ this can be written as
\[\exists R_{l_1},R_{l_2}: ((i_1,j_1) \in R_{l_1} \vee (i_2,j_2) \in R_{l_1}) \wedge((i_1,j_1) \in R_{l_2} \vee (i_2,j_2) \in R_{l_2}).\]
To remove the quantifier, one of the brackets satisfies the equivalence \[\exists R_{l_1}: ((i_1,j_1) \in R_{l_1} \vee (i_2,j_2) \in R_{l_1}) \iff \bigvee_{k=1}^t (x_{i_1,j_1,k} \vee x_{i_2,j_2,k}).\] The same can be done with $l_2$. We just need to make sure, that $l_1 \neq l_2$. We do this, by removing every rectangle one at a time and check the disjunction. In the case, where $l_1$ is removed, there has to be an $l_2$, so that the disjunction is true. We end up with Expression~\eqref{eq4}.\\

\section{Numerical Results}
Firstly, we present some modifications to our models. Although we do not use the MIP formulation for our calculations, we present an additional constraint to better illustrate the changes made to the SAT model. Secondly, we present the matrices, for which we compute the lower bounds. Finally, we present our results.

\subsection{Model Implementation}
In this section, we present the model, that was given to the solver. We modify the versions presented earlier by adding some conditions. The modifications are mostly redundant, but decrease run time. A common approach is to break symmetries~\cite{anders_orbitopal_2026}.
\paragraph{MIP}
Our goal is to cover the ambiguity with differently ordered rectangles. 
We use the lexicographic ordering 
\begin{equation}
     \sum_{i=1}^m 2^{m-i} \cdot u_{i,k^{-}} \geq \sum_{i=1}^m 2^{m-i} \cdot u_{i,k^{-}+1},\quad \forall k^{-} \in [t-1] \label{mip_eq5}
\end{equation}
which is a common way to break symmetries when using binary variables~\cite{kaibel_orbitopal_2011,van_doornmalen_unified_2025}. It was shown by~\cite{friedman_fundamental_2007}, that this indeed yields a strict ordering for binary variables.
This approach gives each rectangle a unique identification number (ID). Each index vector $u_{:,k}$ is assigned an ID by using the vector entries as base two digits. Without loss of generality, the rectangles are ordered by this number from big to small.

\paragraph{SAT Model}
The SAT model is extended by auxiliary variables $a \in \{\texttt{F},\texttt{T}\}^{n\times (t-1)}$ and by five types of inequalities. There are two more constraints using $a$ for breaking symmetry, a redundant one, that worked well in practice, and two, that enforce maximality of rectangles.
\begin{align}
    a_{i^{-},k^{-}} \iff (u_{i^{-},k^{-}} \wedge \lnot u_{i^{-},k^{-}+1})&,\quad \forall (i^{-},k^{-}) \in [m-1]\times[t-1] \label{eq5} \\
    \left(\bigvee_{j=1}^{i-1} a_{j,k^{-}}\right) \vee (u_{i,k^{-}+1} \Rightarrow u_{i,k^{-}}) &,\quad \forall (i,k^{-}) \in [m]\times[t-1]\label{eq6}\\
    v_{j,k} \iff \bigwedge_{i=1}^m (u_{i,k} \Rightarrow (A_{i,j} > 0)) &, \quad \forall (j,k) \in [n]\times [t]\label{eq7}\\
    u_{i,k} \iff \bigwedge_{i=1}^n (v_{j,k} \Rightarrow (A_{i,j} > 0)) &, \quad \forall (i,k) \in [m]\times [t]\label{eq8}\\
    (A_{i_1,j_2} = 0 \vee A_{i_2,j_1} = 0) \implies \lnot(x_{i_1,j_1,k}) \vee \lnot(x_{i_2,j_2,k})
    &, \quad \forall (i_1,i_2,j_1,j_2,k) \in [m]^2\times [n]^2\times [t].\label{eq9}
\end{align}
The symmetry breaking for the SAT model is inspired by~\cite{cbright_harveys_2024}. For the ordering of the rectangles, we aim for a condition equivalent to \[\sum_{i=1}^m 2^{m-i}u_{i,k} \geq \sum_{i=1}^m 2^{m-i}u_{i,k+1}.\] This means, starting from the smallest index $i=1$ the inequality $u_{i,k}\geq u_{i,k+1}$ is evaluated. In case of equality, we look at $i+1$. If it is strict, every subsequent $i$ can be ignored. In our formulation $a_{i^{-}}$ is true, if the inequality at $i^{-}$ is strict (Expression~\ref{eq5}). So any subsequent inequality has to be checked, if no prior $a_{i^{-}}$ is true (Expression~\ref{eq6}).\\
Additionally, we use the Constraint \eqref{eq8}. The implication $v_{j,k} \Rightarrow  \forall i\in [m]: u_{i,k} \Rightarrow (A_{i,j} > 0)$ is redundant. When a certain column $j$ belongs to rectangle $k$, that means all entries are $>0$ if they are also in the rows of $k$. If they are not in row $i$, they can take any non-negative value. The other implication is not redundant. We add a column any time we can. This means, without loss of generality, we only look at maximal rectangles, as was already described in~\cite{gillis_nonnegative_2020}. If there is a solution with non-maximal rectangles, a solution with only maximal rectangles can be obtained by adding rows or columns to non-maximal rectangles.
We use this approach due to its superior empirical runtime performance. In Condition~\ref{eq8} the same is done for rows.\\
Finally, for Constraint~\ref{eq9} we introduce fooling sets.
\begin{definition}
    A fooling set $F\subseteq \mathrm{supp}(S)$ is a set, that satisfies $S_{i_1,j_2} = 0$ or $S_{i_2,j_1} = 0$ for every distinct pair $(i_1,j_1), (i_2,j_2) \in F$. The cardinality of the biggest fooling set is called fooling set number $\omega(S)$.
\end{definition}
The fooling set is a lower bound for RCB$(S)$, because any rectangle contains at most one element of any fooling set. To cover $\omega(S)$ elements in the support, that can not share rectangles, at least $\omega(S)$ rectangles are needed. For the last constraint, we simply take two element fooling sets into account. This means, $(i_1,j_1)$ and $(i_2.j_2)$ can not be in the same rectangle, if $A_{i_1,j_2}$ or $A_{i_2,j_1}$ are zero.\\
\subsection{Non-negative Matrices}
Here, the tested slack matrices are presented. Additionally, there is a non-negative matrix that is not a slack matrix, for which the bounds are computed. We refer to \url{sites.google.com/site/exactnmf/data-set} for the choice of matrices.
\paragraph{Regular Polygons}
Firstly, we define the slack matrix $S_n \in \mathbb{R}_{\geq 0}^{n\times n}$ of the regular $n$-gon with vertices on the unit circle. This corresponds to the definition in~\cite{vandaele_linear_2017}. Let \[c_k \coloneqq \cos\left(\frac{\pi}{n}\right)-\cos\left((2k+1)\frac{\pi}{n}\right)\] for $k \in \{0,\dots, \left\lceil \dfrac{n}{2}\right\rceil -1\}$.
So we get \[S_{i,j} = \min\{c_{j-i},c_{i-j-1}\}\]  with the indices read modulo $n$.
Note, that there are two zeros in every row and column, so $2n$ in total.\\

\paragraph{Linear Euclidean Distance Matrices}
These are non-negative matrices of rank 3 (\cite{gillis_geometric_2012}), but no slack matrices (\cite{gouveia_which_2013}). They are used in~\cite{hrubes_nonnegative_2012} and~\cite{beasley_real_2009} as an example for matrices, which have small real but large non-negative rank. An in-depth discussion of applications can be found in~\cite{anjos_euclidean_2012}. The matrices $M_n \in \mathbb{Z}_{\geq 0}^{n\times n}$ are defined by \[M_{i,j}\coloneqq (a_i-a_j)^2\] for every $i,j \in [n]$ and $a\in \mathbb{R}^n$. We look at the special case $a_i \coloneqq i$. They are only zero on the main diagonal.
\paragraph{Icosidodecahedron}
This is a three dimensional polytope with regular triangles and pentagons as facets. The slack matrix is a $32\times 30$ matrix with 120 zeros.
\paragraph{Correlation Polytope}
The slack matrix $C_n \in \mathbb{R}_{\geq 0}^{2^n\times 2^n}$ to the correlation polytope is indexed by binary vectors $a,b$ of size $n$. The entries are defined as \[C_n(a,b) \coloneqq (1-a^T b)^2.\]
\paragraph{Unique-disjointness matrices}
The unique-disjointness (UDISJ) matrices are similar to $C_n$. They are indexed the same. They differ only in entries, where $a^T b >1$. At those positions, it does not matter if the entry is covered or not. In the model this means, that the corresponding conditions are omitted. Of course, here the RCB is smaller than the one for $C_n$.

\begin{table}[h]
\centering
\setlength{\tabcolsep}{12pt}
\renewcommand{\arraystretch}{1.05}
\caption{Presentation of the runtime of SAT formulation for the rectangle covering bound (RCB). Reported times are in seconds.
Runtimes are split, so SAT refers to finding a feasible rectangle cover of size $\mathrm{RCB}(X)$, and UNSAT denotes the time needed to certify that no rectangle cover of size $\mathrm{RCB}(X)-1$ exists. If the computations did not finish and there is a value for the RCB given, it is known from other bounds (see Section \ref{known_values}) or computations in \cite{baeckelant_computing_2026}.}
\label{tab:rcb_time}
\begin{subtable}[t]{0.48\textwidth}
\centering
\footnotesize
\caption{Linear Euclidean distance matrices}
\label{tab:rcb-ledm}
\begingroup
\fontsize{7.5}{8.5}\selectfont
\begin{tabular}{@{}l cc c@{}}
\toprule
Matrix & \multicolumn{2}{c}{Our SAT model} & RCB \\
\cmidrule(lr){2-3}
 & SAT & UNSAT & \\
\midrule
$M_{5}$   & 1 & 1 & 4 \\
$M_{7}$   & 1 & 1 & 5 \\
$M_{8}$   & 1 & 1 & 5 \\
$M_{10}$  & 1 & 1 & 5 \\
$M_{11}$  & 1 & 5 & 6 \\
$M_{15}$  & 1 & 5 & 6 \\
$M_{20}$  & 1 & 4 & 6 \\
$M_{21}$  & 1 & - & 7 \\
\bottomrule
\end{tabular}
\endgroup
\end{subtable}
\hfill
\begin{subtable}[t]{0.48\textwidth}
\centering
\footnotesize
\caption{Slack matrices of regular $n$-gons}
\label{tab:rcb-ngon}
\begingroup
\fontsize{7.5}{8.5}\selectfont
\begin{tabular}{@{}l cc c@{}}
\toprule
Matrix & \multicolumn{2}{c}{Our SAT model} & RCB \\
\cmidrule(lr){2-3}
 & SAT & UNSAT & \\
\midrule
$S_{5}$  & 1 & 1 & 5 \\
$S_{7}$  & 1 & 1 & 6 \\
$S_{10}$ & 1 & 1 & 7 \\
$S_{12}$ & 1 & 1 & 7 \\
$S_{19}$ & 1 & 1 & 7 \\
$S_{20}$ & 1 & 3 & 7 \\
$S_{21}$ & 1 & 41 & 7 \\
$S_{22}$ & 1 & 118 & 8 \\
$S_{23}$ & 2 & 155 & 8 \\
$S_{25}$ & 2 & 156 & 8 \\
$S_{30}$ & 50 & 174 & 8 \\
$S_{31}$ & 45 & 150 & 8 \\
$S_{32}$ & 676 & 159 & 8 \\
$S_{33}$ & 12390 & 236 & 8 \\
$S_{34}$ & - & 181 & 8 \\
\bottomrule
\end{tabular}
\endgroup
\end{subtable}

\vspace{0.8em}

\begin{subtable}[t]{0.48\textwidth}
\centering
\footnotesize
\caption{Polyhedral slack matrices}
\label{tab:rcb-extra-slack}
\begingroup
\fontsize{7.5}{8.5}\selectfont
\begin{tabular}{@{}l cc c@{}}
\toprule
Matrix & \multicolumn{2}{c}{Our SAT model} & RCB \\
\cmidrule(lr){2-3}
 & SAT & UNSAT & \\
\midrule
Icosidodec. & 358 & $\sim 17$d & $11$ \\
\bottomrule
\end{tabular}
\endgroup
\end{subtable}
\hfill
\begin{subtable}[t]{0.48\textwidth}
\centering
\footnotesize
\caption{Unique-disjointness matrices}
\label{tab:rcb-udisj}
\begingroup
\fontsize{7.5}{8.5}\selectfont
\begin{tabular}{@{}l cc c@{}}
\toprule
Matrix & \multicolumn{2}{c}{Our SAT model} & RCB \\
\cmidrule(lr){2-3}
 & SAT & UNSAT & \\
\midrule
$U_2$ & 1 & 1 & 3 \\
$U_3$ & 1 & 1 & 6 \\
$U_4$ & 1 & 1 & 9 \\
$U_5$ & 31 & - & 18 \\
$U_6$ & 7625 & - & 27 \\

\bottomrule
\end{tabular}
\endgroup
\end{subtable}

\vspace{0.8em}

\begin{subtable}[t]{0.56\textwidth}
\centering
\footnotesize
\caption{Slack matrices of the correlation polytope}
\label{tab:rcb-corr}
\begingroup
\fontsize{7.5}{8.5}\selectfont
\begin{tabular}{@{}l cc c@{}}
\toprule
Matrix & \multicolumn{2}{c}{Our SAT model} & RCB \\
\cmidrule(lr){2-3}
 & SAT & UNSAT & \\
\midrule
$C_2$ & 1 & 1 & 3 \\
$C_3$ & 1 & 1 & 7 \\
$C_4$ & 1 & 1 & 13 \\
$C_5$ & - & - & 25\\
\bottomrule
\end{tabular}
\endgroup
\end{subtable}
\end{table}
\newpage
\subsection{Implementation Details}
For our calculations, we use Kissat which performed well in competitions~\cite{biere_cadical_2024}. The required conjunctive normal form is generated using PySAT. The Kissat instances were run on our server with an AMD EPYC 9654 and 512 GB of RAM.\\

\subsection{Runtime Comparison}

In~\cite{baeckelant_computing_2026} three models have been compared. There is
\begin{enumerate}[label=(\roman*)]
    \item an enumeration approach, which searches for the smallest subset from all maximal rectangles, that is a rectangle covering,
    \item a MILP from~\cite{dewez_computational_2022} that uses two inequalities to ensure the rectangle property for any two elements, and
    \item the MIP approach, we also presented in Section~\ref{sec_mip}.
\end{enumerate}
We present our computed lower bounds together with their runtime. We loosely compare our runtimes to the ones from approach (\romannum{3}) in~\cite{baeckelant_computing_2026}, but do not go into too much detail, because we do not redo their computations on our hardware. We also do not use multiple seeds. Some experiments showed negligible differences of multiple runs on the same model.
In general, we do not impose a time limit, when computing new bounds. The longest runtime was around 17 days for the RCB UNSAT instance of the icosidodecahedron. Any runs without result have been terminated after a similar timespan.

\subsection{RCB}

The runtimes reported in Table~\ref{tab:rcb_time} are given in seconds, with values rounded up, or in days, when to big. If a run did not finish in time, this is indicated by a ``$-$''. Note that the search during SAT time only takes covers of a fixed size into account. We did not bother to look for bigger coverings, because in this case, we knew the size of coverings, we are looking for. Upper bounds are tight in some instances, and the time to find non-optimal coverings empirically tends to be very short.  Because of it's negligible time consumption, generation of the CNF is not part of the runtime either. \\
The value for the icosidodecahedron is new. In case of the polygons slack matrices our runtime is a lot shorter than the reported runtimes of all methods from~\cite{baeckelant_computing_2026}. For the UDISJ and correlation polytope, our method is faster than the model (\romannum{3}).
For the LEDMs, there is a sharp increase in runtime at $n=21$ for the UNSAT instances. Neither the best model (\romannum{3}) from~\cite{baeckelant_computing_2026} nor our model could be solved in time. We ran this for around two weeks. In the other smaller cases our models runtime is similar to the runtime from (\romannum{3}). They are a few seconds faster for smaller $n$. We are faster for the bigger computed $n$, which could indicate a trend. Note, that the difference is small and could come down to faster hardware.

\subsection{RRCB}

\begin{table}[h]
\centering
\setlength{\tabcolsep}{3pt}
\renewcommand{\arraystretch}{1.05}
\caption{Presentation of the runtime of SAT formulation for the refined rectangle covering bound (RRCB). Reported times are in seconds. 
Runtimes are split, so SAT refers to finding a feasible refined rectangle cover of size $\mathrm{RRCB}(X)$, and UNSAT denotes the time needed to certify, that no refined rectangle cover of size $\mathrm{RRCB}(X)-1$ exists.}
\label{tab:rrcb_time}
\begin{subtable}[t]{0.48\textwidth}
\centering
\footnotesize
\caption{Linear Euclidean distance matrices}
\label{tab:rrcb-ledm}
\begingroup
\fontsize{7.5}{8.5}\selectfont
\begin{tabular}{@{}l cc c@{}}
\toprule
Matrix & \multicolumn{2}{c}{Our SAT model} & RRCB \\
\cmidrule(lr){2-3}
 & SAT & UNSAT & \\
\midrule
$M_{9}$   & 1 & 1 & 6 \\
$M_{12}$   & 1 & 18 & 7 \\
$M_{13}$   & 1 & 18 & 7 \\
$M_{14}$  & 90 & 16 & 7 \\
$M_{15}$  & 74 & 21 & 7 \\
$M_{16}$  & 1 & $\sim 6d$ & 8 \\
$M_{17}$  & 1 & $\sim 7d$ & 8 \\
\bottomrule
\end{tabular}
\endgroup
\end{subtable}
\hfill
\begin{subtable}[t]{0.48\textwidth}
\centering
\footnotesize
\caption{Slack matrices of regular $n$-gons}
\label{tab:rrcb-ngon}
\begingroup
\fontsize{7.5}{8.5}\selectfont
\begin{tabular}{@{}l cc c@{}}
\toprule
Matrix & \multicolumn{2}{c}{Our SAT model} & RRCB \\
\cmidrule(lr){2-3}
 & SAT & UNSAT & \\
\midrule
$S_{5}$  & 1 & 1 & 5 \\
$S_{7}$  & 1 & 1 & 6 \\
$S_{10}$ & 1 & 1 & 7 \\
$S_{12}$ & 1 & 1 & 7 \\
$S_{13}$ & 1 & 8 & 8 \\
$S_{16}$ & 46 & 16 & 8 \\
$S_{19}$ & 135 & 2914 & 9 \\
$S_{20}$ & 331 & 3152 & 9 \\
$S_{21}$ & 308 & 4465 & 9 \\
$S_{22}$ & 495 & 5411 & 9 \\
$S_{23}$ & $\sim$1d & 6330 & 9 \\
$S_{24}$ & $\sim$2d & 9497 & 9 \\
$S_{25}$ & - & 14595 & $\geq$9 \\
$S_{30}$ & - & 18977 & $\geq$9 \\
$S_{40}$ & - & $\sim$1d & $\geq$9 \\
\bottomrule
\end{tabular}
\endgroup
\end{subtable}

\vspace{0.8em}

\begin{subtable}[t]{0.48\textwidth}
\centering
\footnotesize
\caption{Slack matrices of the correlation polytope}
\label{tab:rrcb-corr}
\begingroup
\fontsize{7.5}{8.5}\selectfont
\begin{tabular}{@{}l cc c@{}}
\toprule
Matrix & \multicolumn{2}{c}{Our SAT model} & RRCB \\
\cmidrule(lr){2-3}
 & SAT & UNSAT & \\
\midrule
$C_2$ & 1 & 1 & 3 \\
$C_3$ & 1 & 1 & 7 \\
$C_4$ & 147 & 113 & 13 \\
$C_5$ & - & - & -\\
\bottomrule
\end{tabular}
\endgroup
\end{subtable}
\hfill
\end{table}

Similar to the RCB, we report the RRCB runtimes in Table~\ref{tab:rrcb_time}. The UDISJ matrices are not considered here, because they only have 0--1 entries, which means all rectangle coverings are refined. The RRCB of the correlation polytope does not improve upon the RCB. We were not able to compute the RRCB of the icosidodecahedron. We assume this would close the gap to the upper bound of 12 of the $r_+$.
To the best of our knowledge, we are the first to compute the RRCB for LEDMs. For these matrices the runtime for UNSAT is a lot longer than for SAT. While for the regular $n$-gons the RCB and the RRCB take a lot longer to be computed for the SAT instance. In the latter case, the bounds are tight for all computed matrices, so after showing UNSAT for $\mathrm{RRCB}(S_n)-1$ we are done. The SAT instances do not have to be computed. Although we were not able to show unsatisfiability for 9 refined rectangles, we were able to show it for 8, even for some $n\geq 25$. We stopped at $n=40$, as the resulting bounds do not improve upon existing bounds.

\subsection{Known Values of Lower Bounds and the Non-negative Rank}\label{known_values}
\begin{table}[h]
\centering
\scriptsize
\setlength{\tabcolsep}{12pt}
\renewcommand{\arraystretch}{1.08}
\begin{tabular}{c cc cc cc}
\toprule
& \multicolumn{2}{c}{RCB}
& \multicolumn{2}{c}{RRCB}
& \multicolumn{2}{c}{$r_+(S_n)$} \\
\cmidrule(lr){2-3} \cmidrule(l){4-5} \cmidrule(lr){6-7}
$n$ & LB & UB & LB & UB & LB & UB \\
\midrule
$5$ & $5^{\mathrm{FSB}}$ & $5^{\mathrm{triv}}$ & $5^{\mathrm{V17}}$ & $5^{\mathrm{triv}}$ & $5^{\mathrm{FSB}}$ & $5^{\mathrm{triv}}$ \\
$6$ & $5^{\mathrm{V17}}$ & $5^{\mathrm{G16}}$ & $5^{\mathrm{V17}}$ & $5^{\mathrm{V17}}$ & $5^{\mathrm{RCB}}$ & $5^{\mathrm{V17}}$ \\
$7$ & $6^{\mathrm{V17}}$ & $6^{\mathrm{G16}}$ & $6^{\mathrm{V17}}$ & $6^{\mathrm{F12}}$ & $6^{\mathrm{G12}}$ & $6^{\mathrm{F12}}$ \\
$8$ & $6^{\mathrm{V17}}$ & $6^{\mathrm{G16}}$ & $6^{\mathrm{V17}}$ & $6^{\mathrm{V17}}$ & $6^{\mathrm{RCB}}$ & $6^{\mathrm{V17}}$ \\
$9$ & $6^{\mathrm{V17}}$ & $6^{\mathrm{G16}}$ & $7^{\mathrm{V17}}$ & $7^{\mathrm{V17}}$ & $7^{\mathrm{RRCB}}$ & $7^{\mathrm{V17}}$ \\
$10$ & $7^{\mathrm{V17}}$ & $7^{\mathrm{G16}}$ & $7^{\mathrm{V17}}$ & $7^{\mathrm{V17}}$ & $7^{\mathrm{G12}}$ & $7^{\mathrm{V17}}$ \\
$11$ & $7^{\mathrm{V17}}$ & $7^{\mathrm{G16}}$ & $7^{\mathrm{V17}}$ & $7^{\mathrm{V17}}$ & $7^{\mathrm{G12}}$ & $7^{\mathrm{V17}}$ \\
$12$ & $7^{\mathrm{V17}}$ & $7^{\mathrm{G16}}$ & $7^{\mathrm{V17}}$ & $7^{\mathrm{V17}}$ & $7^{\mathrm{G12}}$ & $7^{\mathrm{V17}}$ \\
$13$ & $7^{\mathrm{V17}}$ & $7^{\mathrm{G16}}$ & $8^{\mathrm{V17}}$ & $8^{\mathrm{F12}}$ & $8^{\mathrm{RRCB}}$ & $8^{\mathrm{F12}}$ \\
$14$ & $7^{\mathrm{V17}}$ & $7^{\mathrm{G16}}$ & $8^{\mathrm{D22}}$ & $8^{\mathrm{F12}}$ & $8^{\mathrm{RRCB}}$ & $8^{\mathrm{F12}}$ \\
$15$ & $7^{\mathrm{V17}}$ & $7^{\mathrm{G16}}$ & $8^{*}$ & $8^{\mathrm{F12}}$ & $8^{\mathrm{G12}}$ & $8^{\mathrm{F12}}$ \\
$16$ & $7^{\mathrm{V17}}$ & $7^{\mathrm{G16}}$ & $8^{*}$ & $8^{\mathrm{V17}}$ & $8^{\mathrm{G12}}$ & $8^{\mathrm{V17}}$ \\
$17$ & $7^{\mathrm{V17}}$ & $7^{\mathrm{G16}}$ & $9^{*}$ & $9^{\mathrm{V17}}$ & $\bf{9^{\mathrm{RRCB}}}$ & $9^{\mathrm{V17}}$ \\
$18$ & $7^{\mathrm{V17}}$ & $7^{\mathrm{G16}}$ & $9^{*}$ & $9^{\mathrm{V17}}$ & $\bf{9^{\mathrm{RRCB}}}$ & $9^{\mathrm{V17}}$ \\
$19$ & $7^{\mathrm{V17}}$ & $7^{\mathrm{G16}}$ & $9^{*}$ & $9^{\mathrm{V17}}$ & $\bf{9^{\mathrm{RRCB}}}$ & $9^{\mathrm{V17}}$ \\
$20$ & $7^{\mathrm{V17}}$ & $7^{\mathrm{G16}}$ & $9^{*}$ & $9^{\mathrm{V17}}$ & $\bf{9^{\mathrm{RRCB}}}$ & $9^{\mathrm{V17}}$ \\
$21$ & $7^{\mathrm{V17}}$ & $7^{\mathrm{G16}}$ & $9^{*}$ & $9^{\mathrm{V17}}$ & $9^{\mathrm{G12}}$ & $9^{\mathrm{V17}}$ \\
$22$ & $8^{\mathrm{D22}}$ & $8^{\mathrm{G16}}$ & $9^{*}$ & $9^{\mathrm{V17}}$ & $9^{\mathrm{G12}}$ & $9^{\mathrm{V17}}$ \\
$23$ & $8^{\mathrm{B26}}$ & $8^{\mathrm{G16}}$ & $9^{*}$ & $9^{\mathrm{V17}}$ & $9^{\mathrm{G12}}$ & $9^{\mathrm{V17}}$ \\
$24$ & $8^{\mathrm{V17}}$ & $8^{\mathrm{G16}}$ & $9^{*}$ & $9^{\mathrm{V17}}$ & $9^{\mathrm{G12}}$ & $9^{\mathrm{V17}}$ \\
$25$ & $8^{\mathrm{V17}}$ & $8^{\mathrm{G16}}$ & $9^{*}$ & $10^{\mathrm{F12}}$ & $9^{\mathrm{G12}}$ & $10^{\mathrm{F12}}$ \\
$26$ & $8^{\mathrm{V17}}$ & $8^{\mathrm{G16}}$ & $9^{*}$ & $10^{\mathrm{F12}}$ & $9^{\mathrm{G12}}$ & $10^{\mathrm{F12}}$ \\
$27$ & $8^{\mathrm{V17}}$ & $8^{\mathrm{G16}}$ & $9^{*}$ & $10^{\mathrm{F12}}$ & $9^{\mathrm{G12}}$ & $10^{\mathrm{F12}}$ \\
$28$ & $8^{\mathrm{V17}}$ & $8^{\mathrm{G16}}$ & $9^{*}$ & $10^{\mathrm{F12}}$ & $9^{\mathrm{G12}}$ & $10^{\mathrm{F12}}$ \\
$29$ & $8^{\mathrm{V17}}$ & $8^{\mathrm{G16}}$ & $9^{*}$ & $10^{\mathrm{F12}}$ & $9^{\mathrm{G12}}$ & $10^{\mathrm{F12}}$ \\
$30$ & $8^{\mathrm{V17}}$ & $8^{\mathrm{G16}}$ & $9^{*}$ & $10^{\mathrm{F12}}$ & $9^{\mathrm{G12}}$ & $10^{\mathrm{F12}}$ \\
$31$ & $8^{\mathrm{V17}}$ & $8^{\mathrm{G16}}$ & $9^{*}$ & $10^{\mathrm{F12}}$ & $10^{\mathrm{G12}}$ & $10^{\mathrm{F12}}$ \\
$32$ & $8^{\mathrm{V17}}$ & $8^{\mathrm{G16}}$ & $9^{*}$ & $10^{\mathrm{V17}}$ & $10^{\mathrm{G12}}$ & $10^{\mathrm{V17}}$ \\
$33$ & $8^{\mathrm{V17}}$ & $8^{\mathrm{G16}}$ & $9^{*}$ & $11^{\mathrm{V17}}$ & $10^{\mathrm{G12}}$ & $11^{\mathrm{V17}}$ \\
$34$ & $8^{\mathrm{V17}}$ & $8^{\mathrm{G16}}$ & $9^{*}$ & $11^{\mathrm{V17}}$ & $10^{\mathrm{G12}}$ & $11^{\mathrm{V17}}$ \\
$35$--$40$ & $8^{\mathrm{V17}}$ & $9^{\mathrm{G16}}$ & $9^{*}$ & $11^{\mathrm{V17}}$ & $10^{\mathrm{G12}}$ & $11^{\mathrm{V17}}$ \\
\bottomrule
\end{tabular}
\caption{Summary of the values known for the slack matrices of the regular $n$-gons.\\ 
\emph{Sources:} 
$^{*}$ = this work;
$\mathrm{F12}$~=~\cite{fiorini_extended_2012};
$\mathrm{G12}$~=~\cite{gillis_geometric_2012};
$\mathrm{F13}$~=~\cite{fiorini_combinatorial_2012};
$\mathrm{O14}$~=~\cite{oelze_computing_2014};
$\mathrm{D21}$~=~\cite{dewez_geometric_2021};
$\mathrm{D22}$~=~\cite{dewez_computational_2022};
$\mathrm{V17}$~=~\cite{vandaele_linear_2017};
$\mathrm{G16}$~=~\cite{goucha_ranks_2016};
$\mathrm{RCB}$ means that the lower bound on $r_+(S_n)$ is obtained from the RCB lower bound;
$\mathrm{RRCB}$ means that the lower bound on $r_+(S_n)$ is obtained from the RRCB lower bound;
}
\label{tab:ngon-finite}
\end{table}

We recreate an appendix table from~\cite{baeckelant_computing_2026} for RCB and RRCB. The values for the RRCB lower bound are computed in~\cite{vandaele_linear_2017},~\cite{dewez_computational_2022} and in this work. The upper bound for the RRCB stems, if not computed explicitly, from the upper bound to the non-negative rank. The logarithmic upper bound to $r_+$ is presented in~\cite{fiorini_extended_2012} and improved in~\cite{vandaele_linear_2017}. The lower bound values are mostly provided by the geometric bound from~\cite{gillis_geometric_2012}. For $n=13,14,17,18,19,20$ the RRCB closes the gap to the upper bound. Only for some small $n$, the RCB is among the best lower bounds. In~\cite{vandaele_linear_2017} a lower bound for the RCB is given, which uses the Sperner theorem and a representation of the support by the non-negative factorization. The upper bound given in~\cite{goucha_ranks_2016} uses a fixed row-wise support to reduce the number of possible factorizations.

Though the RRCB was originally presented in~\cite{oelze_computing_2014}, values for small $n$-gons were first published in~\cite{vandaele_linear_2017}, specifically for $n\leq 13$. The case $n=14$ was solved by~\cite{dewez_computational_2022}. We solved the cases $17\leq n\leq 20$, by showing unsatisfiability for $8$ rectangles \footnote{%
The same cases were solved independendtly \href{https://github.com/ajt60gaibb/OpenProblemsInNLA/blob/main/references/holden-nr01-2026-09-13/report.pdf}{here}. The linked report was published, during preparation of this manuscript.%
}. This suffices, because the upper and lower bounds are tight for the computed instances. For $n\geq 25$ we were able to show unsatisfiability for $8$ rectangles too. This matches the lower bound from~\cite{gillis_geometric_2012} for $21 \leq n\leq 30$. The upper bound to the RRCB stems from the upper bound to the $r_+$. The exact $r_+$ values for $n\geq 25$ are mostly unknown. For $n=31,32$ the bound from~\cite{gillis_geometric_2012} is tight.
Note, that we were able to compute the RCB up to $n=33$ exactly. These values were known before from~\cite{goucha_ranks_2016}. The exact value of $n=35$ remains open, but does not deliver new information on the $r_+$, because the best known lower bound is bigger.

\begin{table}[h]
\centering
\footnotesize
\setlength{\tabcolsep}{6pt}
\renewcommand{\arraystretch}{1.08}
\begin{tabular}{c c c cc}
\toprule
& RCB
& RRCB
& \multicolumn{2}{c}{$r_+(M_n)$} \\
\cmidrule(l){2-2} \cmidrule(l){3-3}\cmidrule(lr){4-5}
$n$ & & & LB & UB \\
\midrule
$5$ & $4^{\mathrm{D81}}$ & $5^{*}$ & $5^{\mathrm{G12}}$ & $5^{\mathrm{triv}}$ \\
$6$ & $4^{\mathrm{D81}}$ & $5^{*}$ & $5^{\mathrm{G12}}$ & $5^{\mathrm{H12}}$ \\
$7$ & $5^{\mathrm{D81}}$ & $6^{*}$ & $6^{\mathrm{G12}}$ & $6^{\mathrm{H12}}$ \\
$8$ & $5^{\mathrm{D81}}$ & $6^{*}$ & $6^{\mathrm{G12}}$ & $6^{\mathrm{H12}}$ \\
$9$ & $5^{\mathrm{D81}}$ & $6^{*}$ & $6^{\mathrm{G12}}$ & $7^{\mathrm{H12}}$ \\
$10$ & $5^{\mathrm{D81}}$ & $6^{*}$ & $7^{\mathrm{G12}}$ & $7^{\mathrm{H12}}$ \\
$11$ & $6^{\mathrm{D81}}$ & $7^{*}$ & $7^{\mathrm{B26}}$ & $7^{\mathrm{H12}}$ \\
$12$ & $6^{\mathrm{D81}}$ & $7^{*}$ & $7^{\mathrm{B26}}$ & $7^{\mathrm{H12}}$ \\
$13$ & $6^{\mathrm{D81}}$ & $7^{*}$ & $7^{\mathrm{RRCB}}$ & $8^{\mathrm{H12}}$ \\
$14$ & $6^{\mathrm{D81}}$ & $7^{*}$ & $7^{\mathrm{RRCB}}$ & $8^{\mathrm{H12}}$ \\
$15$ & $6^{\mathrm{D81}}$ & $7^{*}$ & $7^{\mathrm{RRCB}}$ & $8^{\mathrm{H12}}$ \\
$16$ & $6^{\mathrm{D81}}$ & $8^{*}$ & $\bf{8^{\mathrm{RRCB}}}$ & $8^{\mathrm{H12}}$ \\
$17$ & $6^{\mathrm{D81}}$ & $8^{*}$ & $8^{\mathrm{RRCB}}$ & $9^{\mathrm{H12}}$ \\
\bottomrule
\end{tabular}
\caption{Summary of the values known for the LEDMs.}
\label{tab:ledm-finite}
\end{table}

The RCB for the LEDMs is given in~\cite{baeckelant_computing_2026} by \[\mathrm{RCB} = \min \left\{ r: \binom{r}{\lfloor r/2 \rfloor} \geq n\right\}\] using results from~\cite{caen_boolean_1981}. It results from lower bounding by Sperner's Theorem and upper bounding by a construction. The upper bound \[r_+ (M_n) \leq \min_{k\in\mathbb{N}} (\lceil n/2^k\rceil + 2k)\] follows from a result in~\cite{hrubes_nonnegative_2012}.
The best known lower bound for $n\leq 10$ is from~\cite{gillis_geometric_2012}. The gap at $n=9$ could not be closed with the self-scaled lower bound (SSB) from~\cite{fawzi_self-scaled_2016}, which was computed by~\cite{baeckelant_computing_2026}. It is a non-combinatorial bound including the magnitude of entries. The SSB gives the exact $r_+$ for $n=11,12$. The same values are attained by the RRCB. We were able to compute the RRCB for $n\leq 17$, which improved on the best known lower bound for $n \geq 13$ und thus determines xc$(M_{{16}})=8$.

\begin{table}[h]
\centering
\scriptsize
\setlength{\tabcolsep}{8pt}
\renewcommand{\arraystretch}{1.08}
\begin{tabular}{l cc cc cc}
\toprule
& \multicolumn{2}{c}{RCB}
& \multicolumn{2}{c}{RRCB}
& \multicolumn{2}{c}{$r_+$} \\
 \cmidrule(lr){2-3} \cmidrule(lr){4-5} \cmidrule(lr){6-7}
Polytope & LB & UB & LB & UB & LB & UB \\
\midrule
Icosidodecahedron
& $11^{*}$ & $11^*$
& $11^\mathrm{RCB}$ & $12^{\mathrm{B26}}$
& $11^{\mathrm{RCB}}$ & $12^{\mathrm{B26}}$ \\
\bottomrule
\end{tabular}~\caption{Current values of RCB, RRCB and $r_+$ for the icosidodecahedron. For each quantity, LB and UB denote the best currently known lower and upper bounds, respectively. \\
\emph{Sources:} 
$^{*}$ = this work;
$\mathrm{B26}$ = \cite{baeckelant_computing_2026};
$\mathrm{RCB}$ means that the lower bound on $r_+$ is obtained from the RCB lower bound.}
\label{tab:polytope-bounds}
\end{table}

Table~\ref{tab:polytope-bounds} gives an overview over known bounds of the icosidodecahedron. Our calculation of the RCB improved the lower bound by one.
\newpage
\section{Conclusion}
We were able to show, that boolean satisfiability model works well on combinatorial lower bounds of the non-negative rank for certain non-negative matrices. With our model, we can compute new RCB and RRCB values directly. We narrowed the $r_+$ gap of the icosidodecahedron computing the RCB of $11$. For regular $n$-gons, we can compute the RCB up to $n=33$, although these values were known before, through different bounds. The next unknown value of the RCB remains the one for $S_{35}$.
We extended existing computation of the RRCB for regular $n$-gons. This allowed us to determine the non-negative rank \[r_+(S_n) = 9 \ \text{ for } 17\leq n \leq 20.\] The smallest open value for which neither the RRCB nor the $r_+$ is known exactly is $n=25$. However, we were able to improve the lower bound for the RRCB to $9$ for regular $n$-gons with $25\leq n \leq 40$.
We computed the RRCB for LEDMs as well. The values of the RRCB$(M_n)$ equal the best known lower bounds (\cite{fawzi_self-scaled_2016,gillis_geometric_2012}) computed in~\cite{baeckelant_computing_2026} for $n\leq 12$ and improve the best known lower bound for $13\leq n\leq17$. This way, we were able to close the gap to the upper bound in one case and get \[r_+(M_{16}) = 8.\]
\paragraph{Future Work}
Some questions remain unanswered. Understanding the coverings can lead to new constraints or models, which allow computation of new values. More interesting from a theoretical standpoint is the question: Are the RCB and RRCB monotone regarding the $n$ of regular $n$-gons or LEDMs? This seems intuitive, but to the best of our knowledge, this has not been proven yet. It may even be possible to show RRCB = xc for regular $n$-gons:
\begin{conjecture}
    For slack matrices $S_n$ of regular $n$-gons \[\mathrm{r}_+(S_n)=\mathrm{RRCB}(S_n)\] holds.
\end{conjecture}
Our results support this conjecture, which was already posed in~\cite{vandaele_linear_2017}.

\printbibliography

\end{document}